\documentclass[12pt]{article}
\usepackage[final]{epsfig}	
\usepackage{graphics}
\usepackage{amsmath}
\usepackage{amsfonts}
\usepackage{latexsym}
\usepackage{amssymb}

\newtheorem{lemma}{Lemma}[section]

\newtheorem{remark}[lemma]{Remark}

\newtheorem{theorem}{Theorem}

\newtheorem{coro}[lemma]{Corollary}
\newtheorem{problem}{Problem}

\begin{document}
\newcommand{\eps}{{\varepsilon}}
\newcommand{\proofend}{$\Box$\bigskip}
\newcommand{\C}{{\mathbf C}}
\newcommand{\Q}{{\mathbf Q}}
\newcommand{\R}{{\mathbf R}}
\newcommand{\Z}{{\mathbf Z}}

\title
{Tire track geometry: variations on a theme}
\author{Serge Tabachnikov\thanks{Partially supported by an NSF grant}\\
{\it Department of Mathematics, Penn State University}\\
{\it University Park, PA 16802, USA}\\
e-mail: {\it tabachni@math.psu.edu}}
\date{}
\maketitle
\begin{abstract}
We study closed smooth convex plane curves $\Gamma$ enjoying the following property:
a pair of points $x,y$ can traverse $\Gamma$ so that the distances
between $x$ and $y$ along the curve and in the ambient plane do not change; such curves are
called {\it bicycle curves}. Motivation for this study comes from the problem how to
determine the direction of the bicycle motion by the tire tracks of the bicycle wheels;
bicycle curves arise in the (rare) situation when one cannot determine which way the
bicycle went.

We discuss existence and non-existence of bicycle curves, other than circles, in particular,
obtain restrictions on bicycle curves in terms of the ratio of the length of the arc $xy$
to the perimeter length of $\Gamma$, the number and location of their vertices, etc. We also
study polygonal analogs of bicycle curves,  convex equilateral $n$-gons $P$ whose 
$k$-diagonals all have equal lengths. For some values of $n$ and $k$ we prove the rigidity
result that $P$ is a regular polygon, and for some construct flexible
bicycle polygons. 
\end{abstract}

\section{Introduction and outline of results}

The main motivation for what follows comes from the question: ``Which way did the bicycle
go?" A bicycle leaves two tire tracks on the ground, those of the front and
the rear wheels, and the problem is to determine from this pair of curves the direction
of the motion. See \cite{KVW} for a discussion of the problem and, in particular, a
criticism of Sherlock Holmes' approach to it in ``The Priory School" mystery. See also
\cite{DBN,Fi1,Fi2} for various aspects of tire track geometry.

Following \cite{KVW}, we use the next mathematical model for bicycle motion. The
bicycle is represented by an oriented segment of fixed length, say, $L$, whose end points
are the tangency points of the rear and front wheels with the ground. In the process of
motion the end points traverse smooth  plane curves and the segment always remains
tangent to the trajectory of its rear end point.  Let $\gamma(t)$ and $\Gamma(t)$ be
trajectories of the rear and front wheels and let $T(t)$ denote the unit tangent
vector to the curve $\gamma$. Then
\begin{equation}
\label{frontrear}
\Gamma(t)=\gamma(t)+L T(t)
\end{equation}

Following D. Finn \cite{Fi2}, a pair of curves $\gamma$ and $\Gamma$ is
called {\it ambiguous} if they can serve the trajectories of the rear and front bicycle
wheels when traversed in the two opposite directions. Thus one cannot
determine from this pair of tracks which way the bicycle went. If $(\gamma,\Gamma)$ is an
ambiguous pair and a segment of length $2L$ is tangent to $\gamma$ at its midpoint 
then both end points of the segment lie on $\Gamma$, see figure 1.\footnote{Let us
mention a connection to dual billiards: $\Gamma$ is an invariant curve of the dual
billiard map around $\gamma$; see, e.g., \cite{Tab} concerning the dual billiard
problem.}
 An obvious example of
ambiguous curves is a pair of concentric circles whose radii $r$ and $R$ satisfy
$R^2-r^2=L^2$.

\begin{figure}[ht]
\centerline{\epsfbox{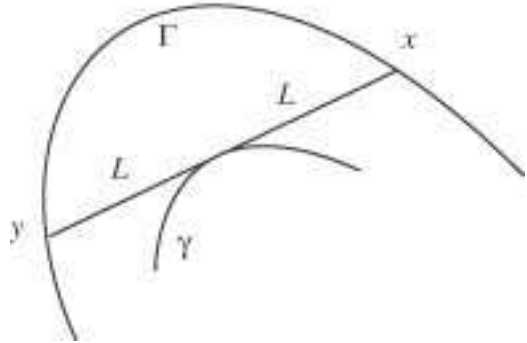}}
\caption{Ambiguous pair of curves}
\end{figure}

The first problem discussed in the present paper is as follows.

\begin{problem}
\label{prblm1}
To describe and study ambiguous pairs of  closed smooth convex plane 
curves $(\gamma,\Gamma)$ where $\Gamma$ is convex. In particular, do
there exist such ambiguous pairs, other than concentric circles?
\end{problem}

One can ask the above question omitting the convexity assumption on $\Gamma$;
it appears, however, that even this restricted version is quite interesting. Problem 1
can be also posed in geometries, other than Euclidean, say, in the hyperbolic plane or
the sphere; one can also ask a similar question in multi-dimensional setting.

Let $(\gamma, \Gamma)$ be an ambiguous pair of curves. We show below in Corollary
\ref{equalall} that, as the segment $xy$ in figure 1 moves around the curve $\gamma$,
the length of the arc $xy$ of the curve $\Gamma$ remains the same. Thus $\Gamma$ has the
following property: a pair of points $x,y$ can traverse the curve so that the distances
between $x$ and $y$ along the curve and in the plane do not change. We will call such
curves {\it bicycle curves}. The ratio of the length of the arc $xy$ to the perimeter
length of $\Gamma$ is called the {\it rotation number} of a bicycle curve and is denoted
by $\rho$. We always assume that the perimeter length of a bicycle curve is $2\pi$.

Our second problem is as follows.

\begin{problem}
\label{prblm2}
 To describe and study closed smooth convex plane bicycle
curves. In particular, what are possible values of the rotation number for non-circular
smooth convex plane bicycle curves?
\end{problem}

Similarly to Problem 1, the questions make sense f in other
geometries and in multi-dimensional setup. 

Problem 2 is more general than Problem 1. Let $\Gamma$ be a
plane bicycle curve. Then the curve $\gamma$ can be recovered as the envelop of the
lines $xy$. However this envelop does not have to be smooth: it may have cusp
singularities, see figure 2. Such singular curves (called wave fronts) still have a well
defined tangent line at every point. In fact, a bicycle can move in such a way that 
the tire track of the rear wheel is a wave front, as in figure 2.  

\begin{figure}[ht]
\centerline{\epsfbox{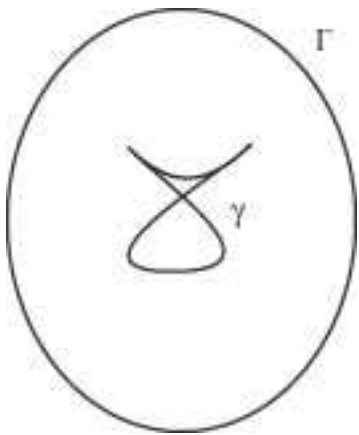}}
\caption{Inner curve is a wave front}
\end{figure}

In higher dimensions, the discrepancy between Problems 1 and 2 gets greater: if
$\Gamma$ is a bicycle curve in, say, 3-dimensional space then the lines $xy$ may not
be tangent to a space curve, and the curve $\gamma$ will not exist at all. 

We also consider polygonal versions of the bicycle curves. Let $P$ be an 
$n$-gon, and let $V_1 V_2 \dots V_n$ be its
consecutive vertices. We understand the indices cyclically, so that $V_{n+1}=V_1$, etc.
A $k$-diagonal is the segment $V_i V_{i+k}$ for some $i=1,\dots ,n$. Let $2\leq k
\leq n/2$. Call $P$ a {\it bicycle $(n,k)$-gon} if it is equilateral, that is, all its 
sides are equal, and all its $k$ diagonals are also equal to each other. An example of a
bicycle 
$(n,k)$-gon is a regular $n$-gon. The ratio $k/n$ will be also called the rotation
number.

The third problem discussed in the paper is as follows.

\begin{problem}
\label{prblm3} 
To describe and study convex plane bicycle
polygons. In particular, for which $n$ and $k$ are there bicycle $(n,k)$-gons,
other than regular $n$-gons?
\end{problem}

Once again, one may ask similar questions for non-convex polygons, for polygons in
multi-dimensional spaces and in geometries, other than Euclidean.

Our results are very far from definitive, and one of the main goals of this paper is to
attract attention to Problems 1--3. 

In Section \ref{examples} we construct examples of
smooth convex plane bicycle curves with the rotation number $\rho=1/2$. There is a
functional space of such curves which, in a sense, are analogous to curves of constant
width. We do not know whether there exists a non-circular smooth convex plane bicycle
curve with the rotation number other than $1/2$. 

Section \ref{vertices} provides restrictions on a smooth convex plane bicycle curve
$\Gamma$ in terms of its vertices (i.e., local maxima or minima of curvature). Let  the
rotation number of $\Gamma$ be $\rho$. Then Theorem \ref{notsmall} states that
every segment of
$\Gamma$ of length $2\pi\rho$ contains a vertex, and Theorem \ref{sixvert} asserts
that the number of vertices of $\Gamma$ is not less than 6. Recall that, by the
celebrated 4-vertex theorem (see, e.g., \cite{Gug}), every simple closed smooth plane
curve has at least 4 vertices. We also give restrictions on the rotation number of a
smooth convex plane bicycle curve $\Gamma$: Theorems \ref{thirdrot} and \ref{quaterrot}
assert that if $\rho=1/3$ or $\rho=1/4$ then $\Gamma$ is a circle. 

Section \ref{modelocking} is devoted to
infinitesimal deformations of the circle in the class of smooth convex plane bicycle
curves. Theorem \ref{deform} describes an interesting mode-locking phenomenon:
 the unit circle has a non-trivial infinitesimal deformation as a
bicycle curve with rotation number $\rho$ if and only if 
\begin{equation}
\label{rotnumrestr}
n\tan (\pi\rho)=\tan(n\pi\rho)
\end{equation}
for some $n\geq 2$. In Section \ref{Defrothalf} we also describe infinitesimal
deformations of the bicycle curves with the rotation number $1/2$ constructed in Section
\ref{examples}.

Section \ref{polygons} concerns bicycle polygons. Theorem \ref{rigidcases}
gives rigidity results: convex plane bicycle $(n,2)$-gons, $(2n+1,3)$-gons, $(2n+1,n)$-gons
and  $(3n,n)$-gons are
regular. On the other hand, Theorem \ref{families} provides a 1-parameter family of
non-regular convex bicycle $(2n,k)$-gons where $k\leq n$ is odd. 

In Section \ref{regular} we describe infinitesimal deformations of regular polygons in the
class of bicycle  polygons. Theorem \ref{regpolyinfdef} asserts that ta regular $n$-gon
admits a non-trivial infinitesimal deformation as a bicycle $(n,k)$-gon if and only if 
$$
\tan\left(kr\frac{\pi}{n}\right) \tan\left(\frac{\pi}{n}\right)
=\tan\left(k\frac{\pi}{n}\right) \tan\left(r\frac{\pi}{n}\right)
$$
for some $2\leq r \leq n-2$. This is an analog of equation (\ref{rotnumrestr}).

Finally, let us mention two papers in which somewhat similar problems are discussed. 

In \cite{Kov}, the following situation is considered. Let $\Gamma$ be a closed convex
plane curve such that two points $x,y$ can traverse $\Gamma$ so that the distance $|xy|$
and the angle $\alpha$ between $xy$ and the tangent line $T_x \Gamma$ remain the same.
The main result of \cite{Kov} is that if 
$\alpha\neq\pi/2$ then $\Gamma$ must be a circle. For $\alpha=\pi/2$, rigidity does not
hold:  $\Gamma$ can be a curve of constant width.

The situation considered in \cite{Gut} is as follows: $\Gamma$ is a smooth
convex closed plane curve such that two points $x,y$ can traverse $\Gamma$ so that the
angles  between $xy$ and the tangent lines $T_x \Gamma$ and $T_y \Gamma$ are both equal to a
constant $\pi\rho$. The main result of \cite{Gut} (see also
\cite{Tab}) is that there exists such a curve, other than a circle, if and only 
equation  (\ref{rotnumrestr}) holds. If the angle is equal to $\pi/2$ then
$\Gamma$ is again a curve of constant width.

Bicycle curves considered here, and the curves studied in \cite{Kov} and
\cite{Gut}, give rise, via duality, to three other classes of smooth convex curves of
interest. The precise notion of duality is the spherical one, namely, the correspondence
``pole--equator" between points of the unit sphere and oriented great circles
(see, e.g., \cite{Arn} or \cite{Tab}). If $x$ and $y$ are two points and $a$ and $b$ are the
corresponding great circles then the spherical distance between $x$ and $y$ equals the
angle between $a$ and $b$. As usual, the duality extends from points and lines to
smooth curves.

\begin{figure}[ht]
\centerline{\epsfbox{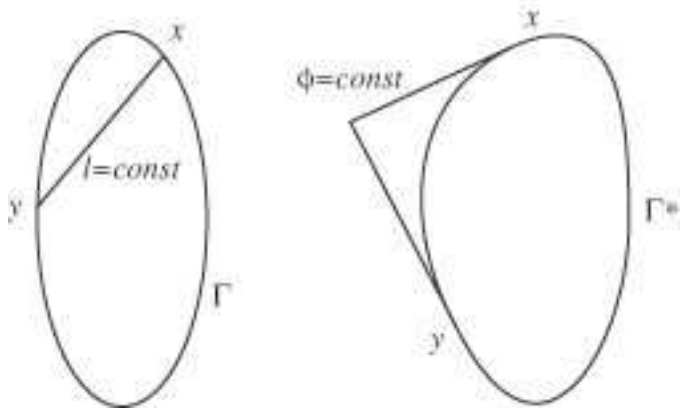}}
\caption{Bicycle curve and its dual curve}
\end{figure}

A curve $\Gamma$, dual to a bicycle curve, has the following property: a pair of 
points $x,y$ can traverse $\Gamma$ so that the distance between $x$ and $y$ along the
curve and the angle between the tangent lines $T_x \Gamma$ and $T_y \Gamma$ do not
change, see figure 3. Likewise, the curves $\Gamma$, dual to the curves studied in
\cite{Gut}, are characterized by the property: a pair of points $x,y$ can traverse
$\Gamma$ so that the distances from points $x$ and $y$ to the intersection point of the
tangent lines $T_x
\Gamma$ and $T_y \Gamma$ remain equal to the same constant. Finally, the  curves
$\Gamma$,  dual to the curves studied in \cite{Kov}, are characterized
by the property: a pair of points $x,y$ can traverse $\Gamma$ so that the distance from
point $x$ to the intersection point of the tangent lines $T_x \Gamma$ and $T_y \Gamma$
and the angle between these tangent lines do not change. A common example for all three
classes is, of course, a circle; are there other examples? 
These three classes of curves can be defined in geometries, other than the spherical
one,  for example, in the Euclidean plane, and their study is an interesting
problem.

{\bf Acknowledgments.} I was introduced to tire track geometry by S. Wagon. 
I have benefited from discussions with my colleagues: M. Bialy, D. Genin, M. Ghomi, E.
Gutkin,  M. Levi, A. Petrunin; I am grateful to them all.

\section{Formul\ae\ for curvature}

Consider two smooth curves $\gamma$ and $\Gamma$ related as in (\ref{frontrear}). 
Let $t$ and $x$ be the arc length parameters on $\gamma$ and $\Gamma$, respectively; 
the correspondence between the two curves is given by 
\begin{equation}
\label{frontrear1}
\Gamma(x(t))=\gamma(t)+L \gamma_t(t).
\end{equation}
Let $k$ and $\kappa$ be the curvatures of
$\gamma$ and $\Gamma$.  Denote by $\alpha(t)$ the angle between the tangent vectors
to the curves at the respective points $\gamma(t)$ and $\Gamma(x(t))$.

\begin{lemma}
\label{curvature}
One has:
$$
-\frac{\pi}{2}<\alpha(t)<\frac{\pi}{2},\qquad \frac{dx}{dt}=\frac{1}{\cos\alpha},
\qquad k=\frac{\tan\alpha}{L}
$$
and
\begin{equation}
\label{curvtract+}
\kappa=\frac{\sin\alpha}{L} + \frac{d\alpha}{dx}.
\end{equation}
\end{lemma}

\noindent {\bf Proof}. Since $t$ is the arc length parameter on $\gamma$, the vector
$\gamma_{tt}$ has magnitude $k$ and is orthogonal to $\gamma_t$. Differentiating
(\ref{frontrear1}), one finds:
\begin{equation}
\label{velocity}
\Gamma_t = \gamma_t + L\gamma_{tt},
\end{equation}
and hence $\Gamma_t \cdot \gamma_t =1$. It follows that $\cos\alpha(t)>0$ for all $t$.
It also follows from (\ref{velocity}) that
\begin{equation}
\label{speed}
\frac{1}{\cos\alpha}=|\Gamma_t|=\sqrt{1+L^2k^2};
\end{equation}
hence $Lk=\tan\alpha$ and $dx/dt=1/\cos\alpha$. 

One has:
\begin{equation}
\label{kappa}
\kappa=\frac{\Gamma_t \times \Gamma_{tt}}{|\Gamma_t|^3}.
\end{equation}
Since $t$ is the arc length parameter on $\gamma$, one has
$
\gamma_{ttt}=(k_t/k)\gamma_{tt}-k^2\gamma_t.
$
Differentiating (\ref{velocity}), one finds:
$
\Gamma_{tt}= \gamma_{tt}+L(k_t/k)\gamma_{tt}-Lk^2\gamma_t.
$
Substitute to (\ref{kappa}) and use
$\gamma_t\times\gamma_{tt}=k$ to obtain:
\begin{equation}
\label{longformula}
\kappa=\frac{k+Lk_t+L^2k^3}{(1+L^2k^2)^{3/2}}=\frac{\sin\alpha}{L}+\alpha_t\cos\alpha
= \frac{\sin\alpha}{L} + \alpha_{x}.
\end{equation}
This completes the proof.
\proofend 

\begin{remark} {\rm As we mentioned in Introduction, one needs to consider a wider class
of curves $\gamma$, namely, wave fronts. The direction of $\gamma$ changes to the
opposite in a cusp, the curvature $k$ at a cusp becomes infinite and changes sign. 
At a cusp point, $\alpha=\pm\pi/2$, and formula (\ref{curvtract+}) still holds.}
\end{remark}

Let $\gamma$ and $\Gamma$ be closed rear and front bicycle tire tracks, that is, closed
smooth curves related by (\ref{frontrear1}). Lemma \ref{curvature} implies that the rear
track is always shorter than the front one.

\begin{coro}
One has:
$$0< {\rm length}\ \Gamma - {\rm length}\ \gamma <L\int_{\gamma} |k|dt.$$
In particular, if $\gamma$ is  convex  then
$${\rm length}\ \Gamma - {\rm length}\ \gamma <2\pi L.$$
\end{coro}

\noindent {\bf Proof}. 
One has:
$$0\leq \sqrt{1+L^2 k^2}-1\leq L|k|.$$
Integrate, using (\ref{speed}), to obtain non-strict inequalities.
Since $|k|>0$ on an open interval, both inequalities are strict.
For a convex curve, $k\geq 0$ and $\int k dt = 2\pi$. 
\proofend

The difference between the lengths of the rear and front tire tracks may be arbitrarily
small: if $\gamma$ is a circle of radius $R$ then 
$$
{\rm length}\ \Gamma - {\rm length}\ \gamma = 2\pi\left(\sqrt{R^2+L^2}-R\right) <
\frac{\pi L^2}{R},
$$
which is small for large $R$.

Now consider the curve 
$
{\bar \Gamma}=\gamma-L\gamma_t,
$
and let $\beta(t)$ be the angle between the tangent vectors to $\gamma$ and ${\bar
\Gamma}$ at the respective points. Up to  orientation reversing,  ${\bar \Gamma}$ is the
trajectory of the front wheel when the rear one traverses $\gamma$ in the opposite
direction, see figure 4.

\begin{figure}[ht]
\centerline{\epsfbox{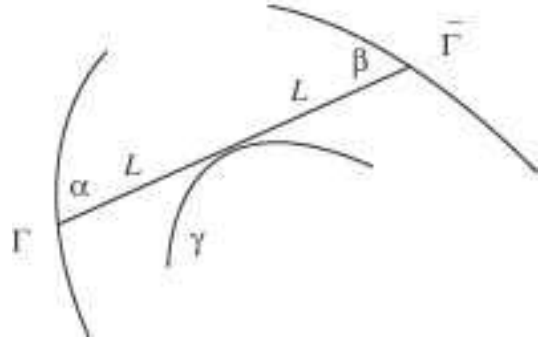}}
\caption{Curves $\Gamma$ and $\bar \Gamma$}
\end{figure}

\begin{lemma}
\label{other}
One has:
$|{\bar \Gamma}_t|=|\Gamma_t|$ and $\beta(t)=\alpha(t).$
The orientation of the frames $({\bar \Gamma}_t, \gamma_t)$ and $(\gamma_t, \Gamma_t)$
coincide.  The curvature of ${\bar \Gamma}$ is given by the formula
\begin{equation}
\label{curvtract-}
\kappa=\frac{\sin\alpha}{L} - \frac{d\alpha}{dx}.
\end{equation}
where $x$ is the arc length parameter on ${\bar \Gamma}$.
\end{lemma}

\noindent {\bf Proof}. 
The computations of Lemma \ref{curvature}, with $L$ replaced by $-L$, yield the formula
for $|{\bar \Gamma}_t|$. Next, ${\bar \Gamma}_t = \gamma_t -
L\gamma_{tt}$, therefore ${\bar \Gamma}_t \cdot\gamma_t=1$, and hence
$\cos\beta=\cos\alpha$.  On the other hand,
$\gamma_t\times{\bar \Gamma}_t=-\gamma_t\times\Gamma_t$, whence the statement concerning
the orientations. The computation of
$\kappa$ proceeds as in Lemma \ref{curvature}, with $\alpha$  replaced by
$-\alpha$, reflecting the orientation, and $L$ by $-L$.
\proofend

Now let $(\gamma, \Gamma)$ be an ambiguous pair of curves; then the curves $\Gamma$ and
${\bar \Gamma}$ coincide. Let $x=\gamma-L\gamma_t$ and $y=\gamma+L\gamma_t$.
One has the following corollary of Lemma \ref{other}.

\begin{coro}
\label{equalall}
The distance between points $x$ and $y$ along the curve $\Gamma$ remains constant and  
the segment $xy$ makes equal angles with  $\Gamma$.
\end{coro}

Let $\Gamma(x)$ be a smooth convex bicycle curve with the rotation number $\rho$ and the
chord of length $2L$, parameterized by arc length. Set $\omega=\pi\rho$. Then the length
of the arc of $\Gamma$ subtended by each chord of length $2L$ equals $2\omega$.\footnote{One
has: $L \leq \sin \omega$ with equality only for a circle. This follows from \cite{Abr}
where it is proved that the average length of the chord subtended by an arc of a fixed
length is not greater than that for a circle.} As before, let $\alpha(x)$ be the angle 
between the segment $\Gamma(x) \Gamma(x+2\omega)$ and the curve $\Gamma$ at points 
$\Gamma(x)$ and $\Gamma(x+2\omega)$ (the angles are equal by Lemma \ref{equalall}).

\begin{theorem}
\label{maineq}
For every $x\in[0,2\pi]$, one has:
\begin{equation}
\label{mastereq}
\sin\alpha(x+\omega)-\sin\alpha(x-\omega)=L\left(\alpha'(x+\omega)+\alpha'(x-\omega)
\right).
\end{equation}
\end{theorem}
\medskip

\noindent {\bf Proof}. 
By Lemmas \ref{other} and \ref{curvature}, the curvatures of $\Gamma$ at points 
$\Gamma(x)$ and $\Gamma(x+2\omega)$
are
\begin{equation}
\label{curvatures}
\kappa(x)=\frac{\sin\alpha(x)}{L} - \alpha'(x),\quad
\kappa(x+2\omega)=\frac{\sin\alpha(x)}{L} +\alpha'(x)
\end{equation} 
where prime denotes $d/dx$. Combining these formulas and shifting $x$ by $\omega$,
yields (\ref{mastereq}).

Alternatively,  compute the total curvature of the arc $\Gamma(x) \Gamma(x+2\omega)$. 
On the one hand, the total turn of this arc is $2\alpha(x)$, on the other hand, due to
(\ref{curvtract-}), it is equal to 
$$
\int_{x}^{x+2\omega} k(\tau) d\tau = \frac{1}{L} \int_{x}^{x+2\omega} \sin \alpha(\tau)
d\tau -\alpha(x+2\omega)+\alpha(x).
$$
Hence
$$
L(\alpha(x+2\omega)+\alpha(x))=\int_{x}^{x+2\omega} \sin \alpha(\tau),
$$
and (\ref{mastereq}) follows by differentiation.
\proofend

By equation (\ref{curvtract+}), a constant solution to (\ref{mastereq}) corresponds to
the unit circle. We are led to the following problem.

\begin{problem}
\label{prblm4}
To describe smooth functions $\alpha(x)$ on the circle $\R/2\pi\Z$
satisfying equation (\ref{mastereq})
\end{problem}

\section{Constructing bicycle curves with the rotation number 1/2}
\label{examples}

Consider an oriented  segment of a fixed length $2L$ in a Euclidean space, characterized
by its midpoint $x$ and the unit vector along the segment $v$. Let $x(t),v(t)$ be smooth
functions describing motion of the segment.

\begin{lemma}
\label{motion}
The endpoints of the segment have equal speeds if and only if $x' \cdot v'=0$. In
the plane case, $x' \cdot v'=0$ if and only if either $v'=0$ or $x'$ is collinear with
$v$.
\end{lemma} 

\noindent {\bf Proof}. 
The endpoints are the vectors $x\pm Lv$ and their velocities are $x'\pm Lv'$. These two
vectors have equal magnitudes if and only if $x' \cdot v'=0$. Since $v$ is a unit vector,
$v'\cdot v=0$. Thus $x'$ and $v$ are perpendicular to $v'$; if $v'\neq 0$ then, in the
plane case,
$x'$ and $v$ are collinear. 
\proofend

\begin{coro}
\label{planemotion}
If a segment of a fixed length is moving in the plane with non-zero angular speed in
such a way that its endpoints have equal speeds then it remains tangent to the
trajectory of its midpoint.
\end{coro}

\noindent {\bf Proof}. 
Since the angular speed does not vanish, $v'\neq 0$. Hence $x'$ and $v$ are collinear,
as claimed.  
\proofend

The envelop $\gamma$ of the moving segments, i.e., the trajectory of its midpoint, may
have singularities, that is, be a wave front. 

We are ready to construct a convex smooth plane closed bicycle curve $\Gamma$ with the
rotation number $1/2$. The chords of $\Gamma$ bisecting the perimeter will have a fixed
length and envelop a wave front $\gamma$. In terms of the bicycle problem, $\Gamma$ is the
tire track of the front and $\gamma$ that of the rear wheel.

\begin{figure}[ht]
\centerline{\epsfbox{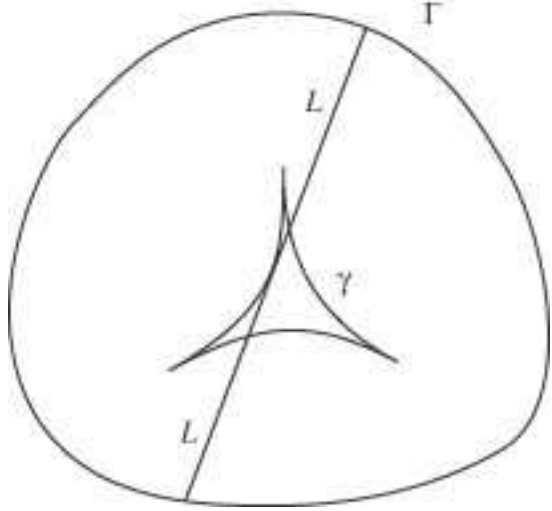}}
\caption{Bicycle curve with the rotation number $1/2$}
\end{figure}

\begin{theorem}
\label{front}
Let $\gamma$ be a closed plane front with an odd number of cusps, total rotation
$\pi$ and without inflections (see figure 5). Let the midpoint of a segment of length
$2L$ traverse $\gamma$ so that the segment remains tangent to $\gamma$. Then, for $L$
large enough,  the endpoints of the segment traverse a convex bicycle curve $\Gamma$. 
\end{theorem} 

\noindent {\bf Proof}. 
When the midpoint of the segment traverses $\gamma$, its endpoints describe one half of
the curve $\Gamma$ each. Since the total rotation of $\gamma$ is $\pi$, these two halves
join smoothly to make a closed curve $\Gamma$.  By Lemma \ref{motion}, the endpoints
move with equal speeds, and hence the segment bisects the perimeter length of $\Gamma$.

To prove that, for $L$ large enough, $\Gamma$ is convex let us, instead of increasing
$L$, rescale $\gamma$ by a small factor  $\lambda$. The curvature of $\gamma$ and its
derivative with respect to the arc length parameter scale as follows:
$
k \mapsto \lambda^{-1} k, k_t \mapsto \lambda^{-2} k_t.
$
By formula (\ref{longformula}), the curvature $\kappa$ of $\Gamma$ scales as follows: 
\begin{equation}
\label{monster}
\frac{k+Lk_t+L^2k^3}{(1+L^2k^2)^{3/2}} \mapsto
\frac{\lambda^{-1}k+\lambda^{-2}Lk_t+\lambda^{-3}L^2k^3}{(1+\lambda^{-2}L^2k^2)^{3/2}}
=\frac{\lambda^{2}k+\lambda Lk_t+L^2k^3}{(\lambda^2+L^2k^2)^{3/2}}.
\end{equation}
Since $\gamma$ has no inflections, $k\neq 0$, and  (\ref{monster}) tends to $1/L$ as 
$\lambda \to 0$. Therefore, for $\lambda$ small enough, (\ref{monster}) is positive.
\proofend

\begin{remark}
{\rm Conversely, it is easy to show that the envelop of the segments that bisect the
perimeter length of a  smooth closed convex plane curve has  total rotation
$\pi$, an odd number of cusps and no inflections.}
\end{remark}

\begin{remark}
{\rm The functions $\alpha(x)$, corresponding to the curves of Theorem \ref{front},
satisfy $\alpha(x-\pi/2)+\alpha(x+\pi/2)=\pi$. Such functions provide solutions to equation
(\ref{mastereq}) with $\omega=\pi/2$ and arbitrary $L$.}
\end{remark}

\section{Restrictions on the number of vertices and on the rotation number}
\label{vertices}

 Assume that $\Gamma$ is a smooth closed convex plane bicycle curve. The
next result shows that if $\Gamma$ is not a circle then the rotation number $\rho$ cannot
be too small. We continue to use the same notation as above, in particular,
$\omega=\pi\rho$.

\begin{theorem}
\label{notsmall}
Every segment of $\Gamma$ of length $2\omega$ contains a vertex.
\end{theorem}

\noindent {\bf Proof}. 
Assume that an arc from $\Gamma(x-\omega)$ to $\Gamma(x+\omega)$ has a non-decreasing and
non-constant curvature. By a Vogt theorem \cite{Gug}, the angle made by the segment
$\Gamma(x-\omega) \Gamma(x+\omega)$ with the curve $\Gamma$ at point $\Gamma(x-\omega)$
is less than that at point $\Gamma(x+\omega)$. This contradicts Lemma \ref{other}.
\proofend

Thus $2\pi\rho$ is not smaller than  the maximal distance between consecutive vertices
along $\Gamma$. Here is another restriction on a bicycle curve $\Gamma$ in terms of
vertices. 

\begin{theorem}
\label{sixvert}
$\Gamma$ has at least 6 vertices.
\end{theorem}

\begin{figure}[ht]
\centerline{\epsfbox{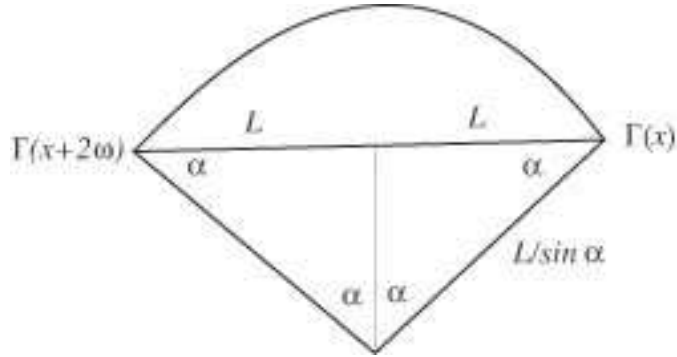}}
\caption{Centers of curvature at $\Gamma(x)$ and $\Gamma(x+2\omega)$ coincide}
\end{figure}

\noindent {\bf Proof}. Let $x$ be a critical point of the function $\alpha$. By formulas
(\ref{curvatures}), the curvatures at points $\Gamma(x)$ and $\Gamma(x+2\omega)$ are
equal. Moreover, it follows from elementary geometry that  the centers of curvature at
these points coincide, see figure 6.  Thus the
circles of curvature at two distinct points of $\Gamma$ coincide. By a Fabricius-Bjerre
theorem \cite{Gug}, $\Gamma$ has at least 6 vertices. 
\proofend

Recall that Theorem \ref{front} provides a variety of bicycle curves with
the rotation number $1/2$. In contrast, we have the next two results.

\begin{theorem}
\label{thirdrot}
If the rotation number equals $1/3$ then $\Gamma$ is a circle.
\end{theorem}

\noindent {\bf Proof}.  For every $x\in [0,2\pi]$, the
triangle $\Gamma(x) \Gamma(x+2\pi/3) \Gamma(x+4\pi/3)$ is equilateral. Let $\alpha(x),
\alpha(x+2\pi/3), \alpha(x+4\pi/3)$ be the respective angles between the chords and the
curve $\Gamma$. The total rotation angle of
$\Gamma$ is
$$
2 (\alpha(x)+\alpha(x+2\pi/3)+\alpha(x+4\pi/3)) = 2\pi
$$
hence
$$
\alpha(x)+\alpha(x+2\pi/3)+\alpha(x+4\pi/3)=\pi.
$$
On the other hand, 
$$\alpha(x)+\alpha(x+2\pi/3)=2\pi/3$$ 
since the angles $\alpha(x)$ and $\alpha(x+2\pi/3)$ sum up to $\pi$ with an interior 
angle of an equilateral triangle. Likewise,
$$
\alpha(x+2\pi/3)+\alpha(x+4\pi/3)=2\pi/3,\quad \alpha(x+4\pi/3)+\alpha(x)=2\pi/3,
$$
therefore $\alpha(x)=\pi/3$ for all $x$. It follows from formula (\ref{curvtract+}) that
the curvature of $\Gamma$ is constant, thus it is a circle.
\proofend

\begin{theorem}
\label{quaterrot}
If the rotation number equals $1/4$ then $\Gamma$ is a circle.
\end{theorem}

\noindent {\bf Proof}. For every $x$, the quadrilateral $\Gamma(x)
\Gamma(x+\pi/2) \Gamma(x+\pi) \Gamma(x+3\pi/2)$ is a rhombus. Let $\alpha(x),
\alpha(x+\pi/2), \alpha(x+\pi), \alpha(x+3\pi/2)$ be the respective angles between the
chords and the curve $\Gamma$. As before, 
$$
2 (\alpha(x)+\alpha(x+\pi/2)+\alpha(x+\pi)+\alpha(x+3\pi/2)) = 2\pi.
$$
Since the opposite angles of a rhombus are equal, one has:
$$
\alpha(x)+\alpha(x+\pi/2)=\alpha(x+\pi)+\alpha(x+3\pi/2)
$$
and
$$
\alpha(x+\pi/2)+\alpha(x+\pi)=\alpha(x+3\pi/2)+\alpha(x).
$$
It follows that 
\begin{equation}
\label{sum}
\alpha(x)+\alpha(x+\pi/2)=\pi/2,
\end{equation}
and hence the rhombus is a square. Therefore the right hand side of equation
(\ref{mastereq}) vanishes, and hence $\sin \alpha(x+\pi/2) = \sin \alpha(x)$. In view of
(\ref{sum}),
$\alpha(x)=\pi/4$, and by formula (\ref{curvtract+}), $\Gamma$ is a circle.
\proofend

\section{Infinitesimal deformations of a circle: mode locking}
\label{modelocking}

Consider a bicycle curve obtained from a circle by an infinitesimal deformation. 

\begin{theorem}
\label{deform}
The unit circle admits a non-trivial infinitesimal deformation as a smooth
closed plane bicycle curve of perimeter $2\pi$ and rotation number $\rho$ if and only if
$\omega$ is a root of the equation
\begin{equation}
\label{modelock}
n \tan (\pi\rho) = \tan(n\pi\rho)
\end{equation}
for some integer $n\geq 2$.
\end{theorem}

\noindent {\bf Proof}. 
Let $\Gamma_0(x)=(\cos x, \sin x)$ be the unit circle and $v(x)$ be a
vector field along $\Gamma_0$. Consider the infinitesimal deformation 
$\Gamma(x)=\Gamma_0(x)+\varepsilon v(x)$, and let it be arc length parameterized as
well. Since $\Gamma'(x)=\Gamma'_0(x)+\varepsilon v'(x)$,  arc length
parameterization is equivalent to  $\Gamma'_0(x) \cdot v'(x) =0$ for all $x$. Thus $(-\sin
x, \cos x) \cdot v'(x) =0$. It follows that $v'(x)= g(x) (\cos x, \sin x)$ for some
function $g$. Note that 
$$\int_0^{2\pi} v'(x) dx =0,$$
hence $g(x)$ is orthogonal to the first harmonics. The differential operator $1+d^2$
``kills" the first harmonics and is a linear isomorphism on the space spanned by all other
harmonics. Hence  $g(x)=f(x)+f''(x)$ for some function $f$. Since $1+d^2$ annihilates the
first harmonics, we assume, without loss of generality, that $f(x)$ is orthogonal to the
first harmonics as well.

One has:
$$
v(x)=\int_0^x (f(\tau)+f''(\tau)) (\cos \tau, \sin \tau) d\tau,
$$
and integration by parts twice gives:
\begin{equation}
\label{fieldv}
v(x)=(f(x) \sin x + f'(x) \cos x  +c_1, -f(x) \cos x + f'(x) \sin x + c_2).
\end{equation}
After a parallel translation, one may set $c_1=c_2=0$. Note also that if $f$ is a
constant then the corresponding deformation $v$ is an infinitesimal rotation of the
circle. Thus, without loss of generality, we assume that $f$ has zero average.

The lengths of the arcs subtended by the chords of $\Gamma$ of constant length
 is $2\omega=2\pi\rho$.  Let this chord length be $L=2(\sin \omega +\varepsilon c)$. Then
$$
|\Gamma(x+\omega)-\Gamma(x-\omega)|=2(\sin \omega +\varepsilon c)
$$ 
which is equivalent to 
$$
(\Gamma_0(x+\omega)-\Gamma_0(x-\omega)) \cdot (v(x+\omega)-v(x-\omega)) = 4c \sin \omega.
$$
A direct computation using (\ref{fieldv}) yields:
$$
(f'(x+\omega)+f'(x-\omega))\sin \omega - (f(x+\omega)-f(x-\omega)) \cos \omega =2c.
$$
The left hand side has zero integral over $[0,2\pi]$, hence $c=0$. One obtains:
\begin{equation}
\label{linear}
(f'(x+\omega)+f'(x-\omega))\sin \omega = (f(x+\omega)-f(x-\omega)) \cos \omega.
\end{equation}
Let
$$
f(\phi)=\sum_{|n|\geq 2} a_n e^{in\phi},\quad {\bar a}_n = a_{-n}
$$
be the Fourier expansion of $f$. Equation (\ref{linear}) is equivalent to
\begin{equation}
\label{Fourier}
a_n (n \cos (n\omega) \sin \omega - \cos \omega \sin(n\omega)) =0
\end{equation}
for all $n\geq 2$.  The Fourier coefficient $a_n$ may be non-zero only if 
$$
n \cos (n\omega) \sin \omega = \cos \omega \sin(n\omega)
$$ 
which is equivalent to (\ref{modelock}).

Conversely, if (\ref{modelock}) holds for some $n\geq 2$, one may choose $f(x)$ to be
a pure $n$-th harmonic, say, $\sin(nx)$. Then the infinitesimal deformation of
the unit circle given by (\ref{fieldv}) is a bicycle curve.
\proofend

\begin{remark}
{\rm Similarly, one may consider infinitesimal deformations of the constant solution
$\alpha(x)=\omega$ with $L=\sin \omega$ of equation (\ref{mastereq}). Such a deformation
exists also if and only if 
\begin{equation}
\label{supereq}
n \tan \omega = \tan (n\omega)
\end{equation} 
for some $n\geq 1$. This does not exclude $n=1$, in which case (\ref{supereq})
trivially holds for all $\omega$. Indeed, a constant solution of equation (\ref{maineq})
admits infinitesimal deformations by the first harmonics; however, such deformations do
not correspond to closed plane curves.}
\end{remark}

\begin{remark}
{\rm 
Note that neither $\rho=1/3$ nor $\rho=1/4$ satisfy (\ref{modelock}), cf.
Theorems \ref{thirdrot} and \ref{quaterrot}. Equation (\ref{modelock}) has solutions for
infinitely many values of $n$. The smallest is $n=4$ for which
$\rho=\arctan(\sqrt{5})/\pi$. }
\end{remark}

\begin{remark}
{\rm 
Note that the roots of equation (\ref{supereq}) are the critical points of the function
$\sin(nx)/\sin x$.}
\end{remark}

\section{Infinitesimal deformations of bicycle curves with rotation number 1/2}
\label{Defrothalf}

Consider a solution of equation (\ref{mastereq}) satisfying
$\alpha(x+\pi)=\pi-\alpha(x)$; such a solution corresponds to a bicycle curve from
Theorem \ref{front} with the rotation number $1/2$. Write
$\alpha(x)=\pi/2+\beta(x)$, then $\beta$ is an odd function.
In this section we study infinitesimal deformations of such solutions.

As the parameter of deformation, we use the change in the rotation number; more precisely,
let $\rho=1/2-\varepsilon/\pi$. Then $\omega=\pi/2-\varepsilon$. Let the deformed
function be $\beta(x)+\varepsilon f(x)$. Without loss of generality, assume that
$f$ is an even function (the odd part of $f$ does not change the rotation number and
can be incorporated into $\beta$). Let the deformed half-length of the chord be
$L+\varepsilon l$.

\begin{theorem}
\label{deformhalf}
There exists an infinitesimal deformation as above if and only if the function $\beta(x)$ 
satisfies the differential equation
\begin{equation}
\label{elliptint}
 L^2 \beta''=(C- \cos \beta)\sin \beta
\end{equation}
where $L$ is as in (\ref{mastereq}), $C$ is a constant, and the respective function $f(x)$
is given by
\begin{equation}
\label{ffrombeta}
Lf=C- \cos \beta.
\end{equation}
\end{theorem}

\noindent {\bf Proof}. Let $x_{\pm}$ stand for $x\pm\pi/2$. Then equation(\ref{mastereq})
can be rewritten as
$$
\cos \left( \beta(x_+-\varepsilon ) + \varepsilon f(x_+-\varepsilon)\right) - 
 \cos \left( \beta(x_-+\varepsilon ) + \varepsilon f(x_-+\varepsilon )\right)=
$$
\begin{equation}
\label{masterrewrite}
(L+\varepsilon l) \left( \beta'(x_+-\varepsilon ) + \varepsilon f'(x_+-\varepsilon
) + \beta'(x_-+\varepsilon ) + \varepsilon
f'(x_-+\varepsilon )\right).
\end{equation}
We compute modulo $\varepsilon^2$. Then
$$
\cos \left( \beta(x_+-\varepsilon ) + \varepsilon f(x_+-\varepsilon)\right)=
\cos \beta(x_+) + \varepsilon (\beta'(x_+)-f(x_+))\ \sin \beta(x_+)
$$ 
and 
$$
\cos \left( \beta(x_-+\varepsilon ) + \varepsilon f(x_-+\varepsilon)\right)=
\cos \beta(x_-) - \varepsilon (\beta'(x_-)+f(x_-))\ \sin \beta(x_-).
$$ 
One also has
$$
\beta'(x_+-\varepsilon )=\beta'(x_+)-\varepsilon \beta''(x_+),\quad 
\beta'(x_-+\varepsilon )=\beta'(x_-)+\varepsilon \beta''(x_-).
$$
Recall that $\beta$ is an odd and $f$ is an even function, hence
$\beta(x_-)=-\beta(x_+)$ and $f(x_-)=f(x_+)$.  Substitute to (\ref{masterrewrite}) and
equate the terms, linear in $\varepsilon$, to obtain:
\begin{equation}
\label{mixedeq}
(\beta'(x)-f(x)) \sin \beta(x) = L(f'(x)-\beta''(x))
\end{equation}
where we replaced $x_+$ simply by $x$.
Equate even and odd parts in (\ref{mixedeq}):
\begin{equation}
\label{evenodd}
 \beta' \sin \beta=Lf', \qquad f \sin \beta= L \beta''.
\end{equation}
The first equation in (\ref{evenodd}) implies: $Lf=C-\cos \beta$. Substitute $f$
to the second equation in (\ref{evenodd}) to obtain the  differential equation
(\ref{elliptint}) on $\beta$.
\proofend

\begin{remark}
{\rm Equation (\ref{elliptint}) can be solved in elliptic integrals; we do not dwell on
this.}
\end{remark}

Theorem \ref{deformhalf} implies that the bicycle curves from
Theorem \ref{front} with the rotation number $1/2$ almost never admit infinitesimal
deformations changing the rotation number. In particular, one has the following
corollary.

\begin{coro}
\label{rigidcircle}
Consider a circle as a bicycle curve with the rotation number $1/2$. Then its every
infinitesimal deformation changing the rotation number yields also a circle.
\end{coro}

\noindent {\bf Proof}.
For a circle, $\beta(x)\equiv 0$. It follows from (\ref{ffrombeta}) that  $f(x)$ is
constant. By (\ref{curvtract+}), the curvature of the deformed curve is also
constant, therefore it is a circle. 
\proofend

\section{Bicycle polygons}
\label{polygons}

We start with a polygonal analog of Lemma \ref{motion}.

\begin{lemma}
\label{motionpoly}
Let $P$ be a plane bicycle $(n,k)$-gon. Then, for every $i$, either the vectors
$V_{i+1}-V_i$ and $V_{i+k+1}-V_{i+k}$ are equal, and then the quadrilateral $V_i V_{i+1}
V_{i+k} V_{i+k+1}$ is a parallelogram, or the quadrilateral
$V_i V_{i+1} V_{i+k} V_{i+k+1}$ is an isosceles trapezoid with the parallel sides
$V_i V_{i+k+1}$ and $V_{i+1} V_{i+k}$. If $P$ is convex then only the latter case is
possible.
\end{lemma}

\noindent {\bf Proof}. The triangles $V_i V_{i+1} V_{i+k}$ and $V_{i+1} V_{i+k} V_{i+k+1}$
are congruent since they have equal corresponding sides. If the segments $V_i V_{i+k}$ and 
$V_{i+1} V_{i+k+1}$ do not intersect then one has the first case of the lemma, and if they
do one has the second case, see figure 7. Clearly only the latter agrees with convexity of
the polygon.
\proofend

\begin{figure}[ht]
\centerline{\epsfbox{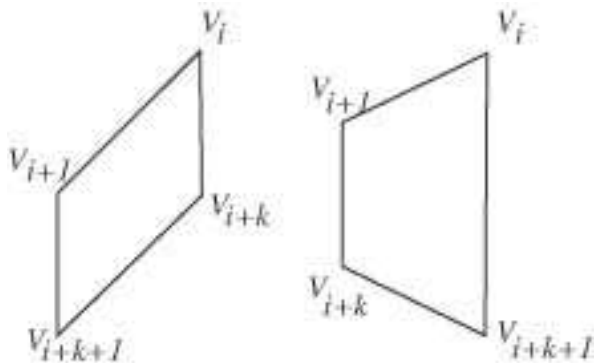}}
\caption{Two cases in Lemma \ref{motionpoly}}
\end{figure}

The next theorem provides some rigidity results on convex bicycle polygons. 

\begin{theorem}
\label{rigidcases} 
In the following cases every convex bicycle $(n,k)$-gon is regular: \\
1) $n$ arbitrary and $k=2$;\\
2) $n$ odd and $k=3$;\\
3) $k$ arbitrary and $n=2k+1$;\\
4) $k$ arbitrary and $n=3k$.
\end{theorem}

\begin{figure}[ht]
\centerline{\epsfbox{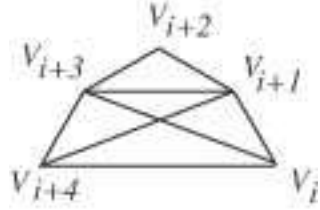}}
\caption{Case of $(n,3)$-gon}
\end{figure}

\noindent {\bf Proof}.
If $k=2$ then the triangles $V_i V_{i+1} V_{i+2}$ are congruent for all $i$, and therefore
all angles of the polygon $P$ are equal. Hence $P$ is regular.

Let $k=3$. Consider 5 consecutive vertices of $P$. Lemma \ref{motionpoly} implies that
the pentagon $V_i V_{i+1} V_{i+2} V_{i+3} V_{i+4}$ has an axis of symmetry passing through
vertex $V_{i+2}$, see figure 8. It follows that the angles at vertices $V_{i+1}$ and
$V_{i+3}$ are equal. This holds for all $i=1,\dots ,n$. If $n$ is odd, this implies that all
angles are equal, and hence $P$ is regular.

Let $n=2k+1$. Then the triangles $V_{i-1} V_i V_{i+k}$ are congruent for all $i$:
two of their sides are $k$-diagonals and the base is a side of the polygon $P$.
Let $\alpha$ be the angle at the base of this triangle. 
The vertices $V_{i-1}, V_i, V_{i+1},V_{i+k}$ and $V_{i+k+1}$ span three such triangles,
and the angle $V_{i-1} V_i V_{i+1}$ is equal to $4\alpha-\pi$, see figure 9. Thus all angles
of $P$ are equal and $P$ is regular. 

\begin{figure}[ht]
\centerline{\epsfbox{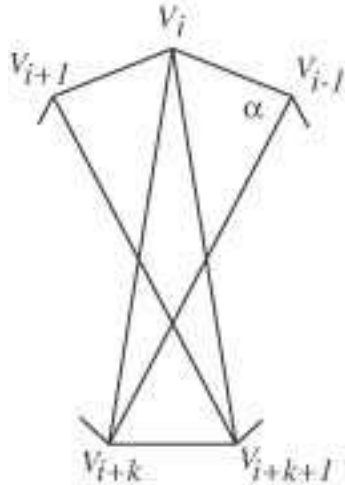}}
\caption{Case of $(2k+1,k)$-gon}
\end{figure}

Finally, consider the case of $n=3k$. We need the following lemma.

\begin{lemma}
\label{twotriangles}
Given two congruent equally oriented
equilateral triangles such that the distances between the
corresponding vertices are equal, one is obtained from
another either by a parallel translation or by a rotation
about its center.
\end{lemma}

\noindent {\bf Proof}. Identify the plane with $\C$ and assume, without loss of generality,
that the first triangle is $(1,q,q^2)$ where $q$ is a cube root of $1$. The motion that
takes 1-st triangle to the 2-nd is given by $z \mapsto uz+v$ where $|u|=1$.
Then one has: 
$$|u-1+v|=|uq-q+v|=|uq^2-q^2+v|,$$
and hence $(1-u)$ is equidistant from the points $v,qv,q^2v$.
Either $v=0$, and the motion is a rotation about the origin, or $1-u=0$, and the motion
is a parallel translation.
\proofend

Now apply Lemma \ref{twotriangles} to a bicycle $3k$-gon $P$. Each triangle $V_i V_{i+k}
V_{i+2k}$ is equilateral. Consider the next triangle $V_{i+1} V_{i+k+1} V_{i+2k+1}$. By
Lemma \ref{twotriangles} and since $P$ is convex, the second triangle is obtained from the
first by a rotation about its center. For all $i$, these rotations have a common center and
equal angles since the sides of $P$ are all equal. It follows that $P$ is a regular polygon.
\proofend

\begin{remark}
{\rm The last case of Theorem \ref{rigidcases} can be viewed as a polygonal analog of Theorem
\ref{thirdrot}.}
\end{remark}

\begin{figure}[ht]
\centerline{\epsfbox{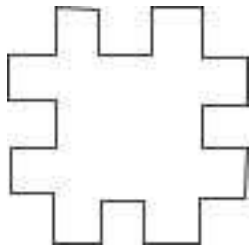}}
\caption{Nonconvex bicycle $(n,2)$-gon}
\end{figure}

\begin{remark}
{\rm If one relaxes the convexity condition on $P$ then rigidity in Theorem \ref{rigidcases}
does not hold anymore. For example, consider a closed $n$-gon on graph paper with unit sides
and all right angles. Clearly there is an abundance of such polygons and they are bicycle
$(n,2)$-gons since their $2$-diagonals have length $\sqrt{2}$, see figure 10. }
\end{remark}

Next we provide some examples of flexible convex bicycle polygons.

\begin{theorem}
\label{families}
For $k$ odd  and $n$ even,  there exists a 1-parameter family of non-congruent
bicycle $(n,k)$-gons.
\end{theorem}

\noindent {\bf Proof}.
Start with a regular $n/2$-gon. Attach to all sides congruent
isosceles triangles to obtain an $n$-gon with equal sides; the altitude of the triangles
is a parameter of the construction. For every odd $k \leq n/2$, the resulting polygon is a
bicycle $(n,k)$-gon  since all $k$-diagonals are congruent by a symmetry
of the original regular $n/2$-gon,  see figure 11 for $(6,3)$ and $(8,3)$-gons.
\proofend

\begin{figure}[ht]
\centerline{\epsfbox{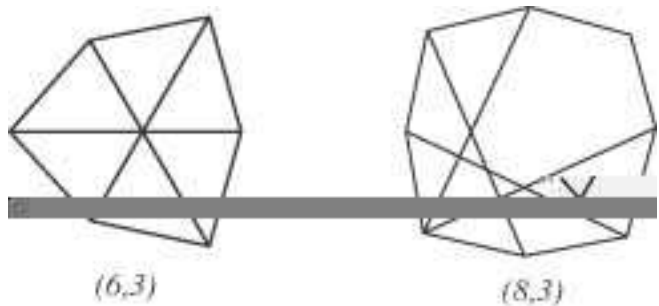}}
\caption{Flexible bicycle polygons}
\end{figure}

\begin{remark}
{\rm Bicycle polygons of Theorem \ref{families} may have the rotation number $1/4$. Thus
there is no polygonal analog of Theorem \ref{quaterrot}.}
\end{remark}

\begin{remark}
{\rm Let $P$ be a convex plane bicycle $(n,k)$-gon. According to Lemma \ref{motionpoly},
$V_i V_{i+1} V_{i+k} V_{i+k+1}$ is an isosceles trapezoid. Consider the circumscribed
circle of this trapezoid and replace the sides $V_i V_{i+1}$ and $V_{i+k} V_{i+k+1}$ by
arcs of this circle, see figure 12. After this is done for all $i$, one obtains a piece-wise
circular bicycle curve  with the rotation number $k/n$. Unless $P$ is regular,
this curve is not differentiable; if $P$ is regular the curve is a circle. This construction
is due to A. Petrunin.}
\end{remark}

\begin{figure}[ht]
\centerline{\epsfbox{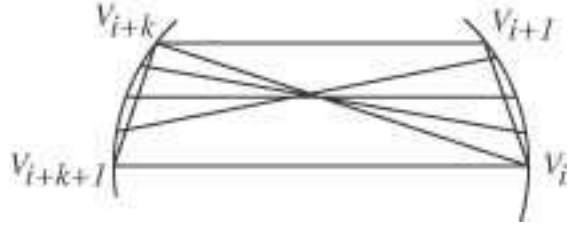}}
\caption{Rotating chord within one trapezoid}
\end{figure}

\begin{remark}
{\rm Up to isometries of the plane, the space of $n$-gons with unit sides is
$n-3$-dimensional. If $k<n/2$ then the condition that all $k$-diagonals are equal
provides $n-1$ relations, and one has an overdetermined system of equations. However, if
$k=n/2$ then there are only $n/2-1$ relations and one has a variety of bicycle
$(2k,k)$-gons.} 
\end{remark}

\section{Infinitesimal deformations of regular polygons}
\label{regular}

This section is a polygonal analog of  Section \ref{modelocking}:
we consider infinitesimal deformations of regular polygons in the class of
bicycle polygons. The main result is as follows.

\begin{theorem}
\label{regpolyinfdef}
A regular $n$-gons admits a
non-trivial infinitesimal deformation as a bicycle $(n,k)$-gon if and only if 
\begin{equation}
\label{polymodecond}
\tan\left(kr\frac{\pi}{n}\right) \tan\left(\frac{\pi}{n}\right)
=\tan\left(k\frac{\pi}{n}\right) \tan\left(r\frac{\pi}{n}\right)
\end{equation}
for some $2\leq r \leq n-2$.
\end{theorem}

\noindent {\bf Proof}.
It will be convenient to denote the angle $\pi/n$ by $\phi$. Consider a regular $n$-gon $P$
whose vertices are
$$
V_i=\left(\cos (2i\phi) ,\ \sin (2i\phi)\right),\quad i=0,\dots ,n-1.
$$
An infinitesimal deformation of $P$ is given by a collection of vectors $U_i$ so that the
vertices of the deformed polygon $P_{\varepsilon}$ are $V_i + \varepsilon U_i$.
As usual, all the computations are modulo $\varepsilon^2$, and the indices are 
understood cyclically.

The polygon $P_{\varepsilon}$ has unit sides:
$$
\left|V_{i+1}-V_i +\varepsilon \left(U_{i+1}-U_i \right)  \right|=1 
$$ 
for all $i=0,\dots ,n-1$, which is equivalent to
$$
\left(V_{i+1}-V_i \right) \cdot \left(U_{i+1}-U_i \right)=0.
$$
Since
$$
V_{i+1}-V_i=2\sin\phi \left(-\sin ((2i+1)\phi) ,\ \cos ((2i+1)\phi)\right)
$$
one has:
$$
U_{i+1}-U_i = t_i \left(\cos ((2i+1)\phi) ,\ \sin ((2i+1)\phi)\right)
$$
for some real $t_i$. Set: 
$$
W_i=\left(\cos ((2i+1)\phi) ,\ \sin ((2i+1)\phi)\right).
$$
Then 
\begin{equation}
\label{UviaW}
U_i=U_0+t_0W_0+t_1 W_1+\dots +t_{i-1} W_{i-1},\quad i=0,\dots ,n-1
\end{equation}
and 
\begin{equation}
\label{consist}
\sum_{i=0}^{n-1} t_i W_i =0.
\end{equation}
Note that adding the same vector to all $U_i$ amounts to parallel translating $P$. Therefore
we may factor out parallel translations by assuming that $U_0=0$.

We can also factor out rotations about the origin. To this end, note that for an
infinitesimal rotation, 
$$
U_i = c \left(-\sin (2i\phi) ,\ \cos (2i\phi)\right)
$$
where $c$ is a constant, and therefore 
$$
U_{i+1}-U_i = C \left(\cos ((2i+1)\phi) ,\ \sin ((2i+1)\phi)\right)
$$
where $C$ is another constant. Hence the rotations correspond to $t_0=t_1=\dots =t_{n-1}$ in
(\ref{UviaW}). To factor the rotations out we assume that 
\begin{equation}
\label{sumt}
\sum_{i=0}^{n-1} t_i =0.
\end{equation}

Next, we consider the condition that all $k$-diagonals of $P_{\varepsilon}$ are equal. As
before, this is equivalent to the equations
\begin{equation}
\label{VdotU}
\left(V_{i+k}-V_i \right) \cdot \left(U_{i+k}-U_i \right)=c
\end{equation}
where $c$ is some constant. A direct computation yields:
$$
W_j \cdot \left(V_{i+k}-V_i \right)= 2\sin(k\phi)\ \sin((2j+1-2i-k)\phi)
$$
and, by (\ref{UviaW}), equation (\ref{VdotU}) can be rewritten as
\begin{equation}
\label{systont}
\sum_{j=0}^{k-1} t_{i+j}\ \sin((2j+1-k)\phi)=C
\end{equation}
where $C$ is another constant. The sum of the left hand sides of (\ref{systont}) over
$i=0,\dots ,n-1$ is zero, therefore $C=0$. One finally obtains the system of linear equations
on the variables $t_i$:
\begin{equation}
\label{systont1}
\sum_{j=0}^{k-1} t_{i+j}\ \sin((2j+1-k)\phi)=0,\quad i=0,\dots ,n-1,
\end{equation}
along with (\ref{consist}) and (\ref{sumt}). The matrix $A$ of system (\ref{systont1}) is
$$
\begin{array}{cccc}
a_0&a_1&\dots&a_{n-1}\\
a_{n-1}&a_0&\dots&a_{n-2}\\
\vdots&\vdots&\ddots&\vdots\\
a_1&a_2&\dots&a_0
\end{array}
$$
where 
$a_j=\sin((2j+1-k)\phi)$ for $j=0,\dots ,k-1$ and $a_j=0$ otherwise.

To study this system let $\xi=\exp(2\phi\sqrt{-1})$ be $n$-th primitive root of unity, and
set:
$$
\theta_r=\sum_{j=0}^{n-1} a_j \xi^{jr},\quad r=0,\dots ,n-1.
$$
Let $B$ be the matrix $b_{ij}=\xi^{i(j-1)}$. Note that $B$ is non-degenerate. The next lemma
is verified by a direct computation that we omit.

\begin{lemma}
\label{circular}
One has:
$$
A=B^{-1}D B.
$$
where
$$
D={\rm Diag} (\theta_{n-1}, \theta_{n-2},\dots ,\theta_0).
$$
\end{lemma}

Next, we need the following result.

\begin{lemma}
\label{thetacomp}
One has: 
\begin{equation}
\label{thetaequals}
2\theta_r= \left(\frac{\sin (k(r+1)\phi)}{\sin ((r+1)\phi)} 
- \frac{\sin (k(r-1)\phi)}{\sin ((r-1)\phi)}\right)
\exp{\left( \sqrt{-1}\left( \frac{\pi}{2}-(k-1)r\phi\right)\right)}
\end{equation}
where, for $r=1$ and $r=n-1$, one has:
$$
\frac{\sin (k(r-1)\phi)}{\sin ((r-1)\phi)}=k.
$$
\end{lemma}

\noindent {\bf Proof}.
The real and imaginary parts of $2\theta_r$ are the following sums:
\begin{equation}
\label{realpart}
\sum_{j=0}^{k-1} \sin ((2j+1-k-2rj)\phi) + \sum_{j=0}^{k-1}\sin ((2j+1-k+2rj)\phi)
\end{equation}
and
\begin{equation}
\label{impart}
\sum_{j=0}^{k-1} \cos ((2j+1-k-2rj)\phi) - \sum_{j=0}^{k-1}\cos ((2j+1-k+2rj)\phi).
\end{equation}
Let us use the identities:
\begin{equation}
\label{sinidentity}
\sum_{j=0}^{k-1} \sin (\alpha+2j\beta)=\frac{\sin (k\beta)\ 
\sin(\alpha+(k-1)\beta)}{\sin\beta}
\end{equation}
and
\begin{equation}
\label{cosidentity}
\sum_{j=0}^{k-1} \cos (\alpha+2j\beta)=\frac{\sin (k\beta)\ 
\cos(\alpha+(k-1)\beta)}{\sin\beta}
\end{equation}
where (\ref{sinidentity}) and (\ref{cosidentity}) are equal to $k\sin\alpha$ and
$k\cos\alpha$ respectively if $\sin\beta=0$.

Apply (\ref{sinidentity}) and (\ref{cosidentity}) to (\ref{realpart}) and (\ref{impart})
with an appropriate choice of $\alpha$ and $\beta$ to to find that (\ref{realpart}) equals 
$$
\left(\frac{\sin (k(r+1)\phi)}{\sin ((r+1)\phi)} 
- \frac{\sin (k(r-1)\phi)}{\sin ((r-1)\phi)}\right) \sin((k-1)r\phi)
$$
and (\ref{impart}) equals
$$
\left(\frac{\sin (k(r+1)\phi)}{\sin ((r+1)\phi)} 
- \frac{\sin (k(r-1)\phi)}{\sin ((r-1)\phi)}\right) \cos((k-1)r\phi).
$$
This is equivalent to the statement of the lemma.
\proofend

Now we can complete the proof of the theorem. We are interested in the system of
linear equations $(B^{-1}D B) \bar t=0$ where $\bar t = (t_0, \dots ,t_{n-1}) \in \C^n$
satisfies (\ref{consist}) and (\ref{sumt}). These two
conditions are equivalent to
$$
\bar t \cdot (1,\xi,\xi^2,\dots ,\xi^{n-1})=0,\quad \bar t \cdot (1,1,\dots ,1)=0,
$$
that is, to the condition that the first and the last components
of the vector $B \bar t$ vanish. Therefore the dimension of the space of solutions of our
system equals the number of zeros among the numbers $\theta_1, \theta_2,\dots ,\theta_{n-2}$.

By Lemma \ref{thetacomp}, this equals to the number of $r=1,2,\dots ,n-2$ such that
\begin{equation}
\label{proportion}
\frac{\sin (k(r+1)\phi)}{\sin ((r+1)\phi)} 
= \frac{\sin (k(r-1)\phi)}{\sin ((r-1)\phi)}.
\end{equation}
If $r=1$ then (\ref{proportion}) becomes
$
\sin (2k\phi)=k\sin (2\phi);
$ 
according to Lemma 6 in \cite{Abr} this has no solutions. By elementary trigonometry, for
$r=2,\dots ,n-2$, equation (\ref{proportion}) can be rewritten as (\ref{polymodecond}), and
this completes the proof.
\proofend

\begin{remark}
{\rm Recall that Theorem \ref{families} provides 1-parameter families of bicycle
$(n,k)$-gons for $n$ even and $k$ odd. For such values of $n$ and $k$, equation
(\ref{polymodecond}) holds for $r=n/2$: indeed, $r\phi=\pi/2$ and $\tan (r\phi) = \tan
(kr\phi) = \infty$. Likewise, if $n=2k$, equation (\ref{polymodecond}) holds for every
odd $r$. Note a curious duality between $k$ and $r$ in (\ref{polymodecond}).}
\end{remark}

\end{document}